\documentclass{article}%
\usepackage{amssymb}
\usepackage{amsmath}
\usepackage{graphicx}
\usepackage{amsfonts}%
\providecommand{\U}[1]{\protect\rule{.1in}{.1in}}
\begin{document}

\title{\textsc{Four Ways from Universal to Particular: }\\How Chomsky's Language-Acquisition Faculty is Not Selectionist}
\author{David Ellerman\\Visiting Scholar\\Philosophy Department\\University of California at Riverside}
\maketitle

\begin{abstract}
Following the development of the selectionist theory of the immune system,
there was an attempt to characterize many biological mechanisms as being
"selectionist" as juxtaposed to "instructionist." But this broad definition
would group Darwinian evolution, the immune system, embryonic development, and
Chomsky's language-acquisition mechanism as all being "selectionist." Yet
Chomsky's mechanism (and embryonic development) are significantly different
from the selectionist mechanisms of biological evolution or the immune system.
Surprisingly, there is an abstract way using two dual mathematical logics to
make the distinction between genuinely selectionist mechanisms and what are
better called "generative" mechanisms. This note outlines that distinction.

\end{abstract}
\tableofcontents

\section{Introduction}

There is a simple, abstract, and logical way to classify four ways that
"universals" can give rise to "particulars." One might expect such an abstract
classification to have little power to make interesting or significant
distinctions. Yet it does give a clear differentiation between what would be
called a \textit{selectionist} mechanism and the \textit{generative} mechanism
of Noam Chomsky's language-acquisition faculty or universal grammar (UG). The
purpose of this note is to spell out that distinction.

\section{Selectionist mechanisms}

There is a long tradition, growing out of biological thought, to juxtapose
"selectionist" mechanisms to "instructionist" (or Lamarckian) mechanisms
(\cite{med:fut}, \cite{jerne:al}, \cite{edel:nd}, \cite{cziko:wm},
\cite{dennett:darwin}). Originally the distinction was drawn in rather general
terms. In an instructionist or Lamarckian mechanism, the environment would
transmit detailed instructions about a certain adaptation to an organism,
while in a selectionist mechanism, a diverse variety of (perhaps random)
variations would be generated, and then some adaptations would be selected by
the environment but without detailed instructions from the environment.

The selectionist-instructionist juxtaposition gained importance with the
development of the selectionist theory of the immune system. \cite{jerne:nst}
There is some "generator of diversity" that generates a wide variety of
possible adaptations, and then interaction with the environment differentially
amplifies some possibilities while the others languish, atrophy, or die off.
For instance, in the case of the human immune system, "It is estimated that
even in the absence of antigen stimulation a human makes at least $10^{15}$
different antibody molecules--its \textit{preimmune antibody repertoire}."
(\cite[p. 1221]{alberts:mbcell} quoted in \cite[p. 187]{jenkins:biolin})
Gerald Edelman has sharpened the selectionist definition and generalized its application.

\begin{quotation}
The long trail from antibodies to conscious brain events has reinforced my
conviction that evolution, immunology, embryology, and neurobiology are all
sciences of recognition whose mechanics follow selectional principles. ...All
selectional systems follow three principles. There must be a generator of
diversity, a polling process across the diverse repertoires that ensue, and a
means of differential amplification of the selected variants. \cite[p.
7367]{edel:recog}
\end{quotation}

\noindent In particular, Edelman develops a selectionist theory of brain development.

\begin{quotation}
[T]he theoretical principle I shall elaborate here is that the origin of
categories in higher brain function is somatic selection among huge numbers of
variants of neural circuits contained in networks created epigenetically in
each individual during its development; this selection results in differential
amplification of populations of synapses in the selected variants. In other
words, I shall take the view that the brain is a selective system more akin in
its workings to evolution than to computation or information processing.
\cite[p. 25]{edel:nd}
\end{quotation}

\noindent The key point is that the possibilities must be in some sense
\textit{actualized} or \textit{realized} (e.g., as antibodies in low
concentration in the immune system) in order for selection to operate on and
differentially amplify or select some of the actual variants while the others
languish, atrophy, or die off.

\section{Selectionist vs. generative approaches to universal grammar}

What would the child's language-learning faculty look like if it was
selectionist in this sense? There is a na\"{\i}ve and a sophisticated
selectionist account of language acquisition. In the na\"{\i}ve account, the
child would (perhaps randomly) generate a diverse range of babblings, some of
which would be differentially reinforced by the linguistic environment
\cite{skinner:aboutbehav}. 

\begin{quotation}
\noindent Skinner, for example, was very explicit about it. He pointed out,
and he was right, that the logic of radical behaviorism was about the same as
the logic of a pure form of selectionism that no serious biologist could pay
attention to, but which is [a form of] popular biology -- selection takes any
path. And parts of it get put in behaviorist terms: the right paths get
reinforced and extended, and so on. It's like a sixth grade version of the
theory of evolution. It can't possibly be right. But he was correct in
pointing out that the logic of behaviorism is like that [of na\"{\i}ve
adaptationism], as did Quine.\cite[p. 53]{chom-mcgil} 
\end{quotation}

\noindent As noted, Willard Van Orman Quine adopted essentially this approach
to language learning.

\begin{quotation}
\noindent An oddity of our garrulous species is the babbling period of late
infancy. This random vocal behavior affords parents continual opportunities
for reinforcing such chance utterances as they see fit; and so the rudiments
of speech are handed down. \cite[p. 73]{quine:wando-new}

It remains clear in any event that the child's early learning of a verbal
response depends on society's reinforcement of the response in association
with the stimulations that merit the response, from society's point of view,
and society's discouragement of it otherwise. \cite[p. 75]{quine:wando-new}
\end{quotation}

\noindent Since language users can generate a variety of rule-based
grammatical sentences never before spoken, the child would have to rather
miraculously generalize the grammatical rules from the reinforced variants. In
order to be adequate as an explanation, a model needs to take seriously the
speaker's rule-based competency.

There is a sophisticated version of a selectionist model for the
language-acquisition faculty or universal grammar (UG) which could be called
the \textit{format-selection (FS) approach} (Chomsky, private communication).
The diverse variants that are actualized in the mental mechanism are different
sets of rules or grammars. Then given some linguistic input from the
linguistic environment, the grammars are evaluated according to some
evaluation metric, and the best rules are selected.

\begin{quotation}
\noindent Universal grammar, in turn, contains a rule system that generates a
set (or a search space) of grammars, $\{G_{1},G_{2},\ldots,G_{n}\}$. These
grammars can be constructed by the language learner as potential candidates
for the grammar that needs to be learned. The learner cannot end up with a
grammar that is not part of this search space. In this sense, UG contains the
possibility to learn all human languages (and many more). ... The learner has
a mechanism to evaluate input sentences and to choose one of the candidate
grammars that are contained in his search space. \cite[p. 292]%
{nowak:format-sel}
\end{quotation}

\noindent The idea is that after a sufficient stream of linguistic inputs, the
mechanism would converge to the best grammar that matches the linguistic
environment. Since it is optimizing over sets of rules, this model at least
takes seriously the need to account for rule-based competency. Early work
(through the 1970s) on accounting for the language-acquisition faculty or
universal grammar (UG) seems to have assumed such an approach.

\begin{quotation}
The earliest ideas were roughly as follows. Suppose that UG provides a certain
format for languages, that is, a specification of permitted types of rules and
permissible interactions among them. Any rule system satisfying the proposed
format qualifies as a possible human language. ... The mind employs certain
primitive operations to interpret some of the data presented to it as
linguistic experience, then selects among the languages consistent with this
experience in accordance with an evaluation metric that assigns an abstract
value to each language. \cite[p. 52]{chom:knowoflang}
\end{quotation}

The problems that eventually arose with the FS approach could be seen as the
conflict between descriptive and explanatory adequacy. In order to describe
the enormous range of human language grammars, the range of grammars
considered would make for an unfeasible computational load of evaluating the
linguistic experience. If the range was restricted to make computation more
feasible, then it would not explain the variety of human languages.

\begin{quotation}
\noindent It was an intuitively obvious way to conceive of acquisition at the
time for -- among other things -- it did appear to yield answers and was at
least more computationally tractable than what was offered in structural
linguistics, where the alternatives found in structural linguistics could not
even explain how that child managed to get anything like a morpheme out of
data. But the space of choices remained far too large; the approach was
theoretically implementable, but completely unfeasible. \cite[p.
173]{chom-mcgil}
\end{quotation}

Instead of the format-selection (FS) approach, the alternative
\textit{principles and parameters (P\&P) approach} (\cite{chom:govt},
\cite{chomlas:pandp}, \cite{chom:minprog}) to universal grammar was then developed:

\begin{quotation}
\noindent we no longer consider UG as providing a format for rule systems and
an evaluation metric. Rather, UG consists of various subsystems of principles;
it has the modular structure that we regularly discover in investigation of
cognitive systems. Many of these principles are associated with parameters
that must be fixed by experience. The parameters must have the property that
they can be fixed by quite simple evidence, because this is what is available
to the child; the value of the head parameter, for example, can be determined
from such sentences as \textit{John saw Bill} (versus \textit{John Bill saw}).
Once the values of the parameters are set, the whole system is operative.
\cite[p. 146]{chom:knowoflang}
\end{quotation}

Our purpose here is to give an abstract conceptual differentiation of the P\&P
approach from the sophisticated selectionist approach of the FS system (not to
mention from the crude selectionism in behaviorism or na\"{\i}ve Darwinism). 

Jerry Fodor and Massimo Piattelli-Palmarini took a different approach to that
differentiation. In the process, they generated some controversy with
evolutionary biologists by claiming that "Skinner's account of learning and
Darwin's account of evolution are identical in all but name" or, to be more
precise, "what is wrong with Darwin's account of the evolution of phenotypes
is very closely analogous to what is wrong with Skinner's account of the
acquisition of learned behavior." \cite[p. xvi]{fodorp-p:darwinwrong} They
emphasize aspects of what is broadly called the "neo-neo-Darwinism," Evo Devo,
or the "extended evolutionary synthesis" \cite{pigl-muller:evo-devo}. Instead
of wading into that controversy, we take the different approach in this note
of showing how selectionist mechanisms (e.g., FS) and "generative" P\&P-type
mechanisms can be differentiated at a very abstract logico-mathematical level.
Hence we must turn to a recent development in mathematical logic.

\section{The two dual forms of mathematical logic}

George Boole \cite{boole:lot} originally developed what might be called
\textit{Boolean logic} as the logic of subsets, not the logic of propositions.
The interpretation solely in terms of propositions and the name
\textit{propositional logic} came later.

\begin{quotation}
\textit{The algebra of logic} has its beginning in 1847, in the publications
of Boole and De Morgan. This concerned itself at first with an algebra or
calculus of classes, to which a similar algebra of relations was later added.
Though it was foreshadowed in Boole's treatment of "Secondary Propositions," a
true propositional calculus perhaps first appeared from this point of view in
the work of Hugh MacColl, beginning in 1877. \cite[pp. 155-156]{church:ml}
\end{quotation}

\noindent When Boolean logic is interpreted as the logic of subsets, then
variables stand for subsets of some given universe set $U$, the operations are
subset operations, and a (\textit{subset-})\textit{valid} formula or
(\textit{subset-})\textit{tautology} is a formula so that no matter what
subsets of $U$ are substituted for the variables, the whole formula will
evaluate to $U$ for any non-empty $U$. It is then a theorem (known to Boole),
not a definition, that it suffices to consider on the case where $U=1$ is a
singleton which has only two subsets $1$ and $\emptyset$ (the empty set).
Hence validity in the special case of propositions with the two truth-values
$1$ and $0$, i.e., \textit{truth-table validity}, is equivalent to general
subset-validity. Eventually, the special case of propositional variables came
to dominate so truth-table validity became the \textit{definition} of a
tautology rather than a theorem about subset-validity (see \textit{any}
contemporary logic textbook).

What is lost by this focus on the special case of propositional logic rather
than the general case of subset logic? Around the middle of the twentieth
century, the theory of categories was formalized \cite{eilmac:gte} and an
older informal notion of duality in algebra was formalized as the
reverse-the-arrows duality of category theory \cite{mac:cwm}. The older
informal duality in algebra was the juxtaposition of subgroups to quotient
groups, subrings to quotient rings, and in general subobjects to quotient
objects--which in the basic case of sets was the juxtaposition of subsets to
quotient sets (the latter being equivalent to equivalence relations or
partitions on a set). For instance, F. William Lawvere calls the general
notion of a subobject a "part" and \textquotedblleft The dual notion (obtained
by reversing the arrows) of `part' is the notion of \textit{partition}%
.\textquotedblright\ \cite[p. 85]{law:sfm} Hence when the special case of
"propositional" logic is seen as the general logic of subsets, then the idea
arises of there being a \textit{dual logic of quotient sets or partitions}
(\cite{ell:partitions}, \cite{ell:intropartlogic}). That idea of a dual logic
does not arise when subset logic is seen only as propositional logic since
"propositions" do not have a category-theoretic dual.

\section{The two lattices of subsets and partitions}

The two logics of subsets and partitions are represented algebraically by the
Boolean algebra of subsets of a universe $U$ and the algebra of partitions on
a universe set $U$ ($\left\vert U\right\vert \geq2$). For our purposes, it
suffices to consider the two lattices, the familiar Boolean lattice of subsets
of $U$ (where the partial order is inclusion) and the lattice of partitions on
$U$ where the partial order is the "refinement" relation between partitions.

A \textit{partition} $\pi=\left\{  B\right\}  $ on a universe $U$ is a set of
nonempty blocks $B$ that are disjoint and whose union is $U$. Given two
partitions $\pi=\left\{  B\right\}  $ and $\sigma=\left\{  C\right\}  $ on
$U$, the partition $\pi$ \textit{refines} the partition $\sigma$, written
$\sigma\preceq\pi$, if for every block $B\in\pi$, there is a block $C\in
\sigma$ such that $B\subseteq C$. Figure $1$ illustrates the two lattices for
the universe $U=\left\{  a,b,c\right\}  $ (where the partial order is
indicated by the lines).%

\begin{center}
\includegraphics[
height=1.8821in,
width=3.8232in
]%
{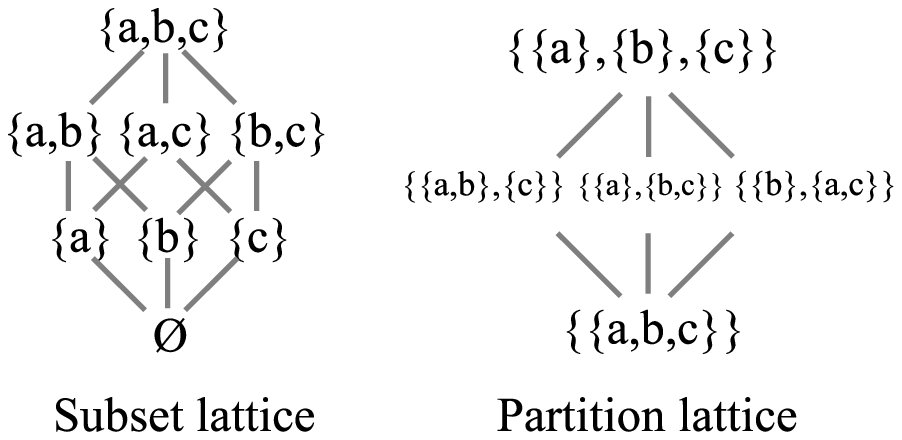}%
\\
Figure $1$: Lattices of subsets and partitions
\end{center}

Each lattice has a top and bottom and that gives us the four universals in our
title. In the Boolean lattice of subsets, the top is the universe set $U$ and
the bottom is the null set $\emptyset$. In the lattice of partitions, the top
is the \textit{discrete partition} $\mathbf{1}=\left\{  \left\{  u\right\}
:u\in U\right\}  $ where all blocks are singletons, and the bottom is the
\textit{indiscrete partition} $\mathbf{0}=\left\{  \left\{  U\right\}
\right\}  $ with only one block consisting of the universe $U$.

The duality between subsets of a set and partitions on a set extends to the
extensive analogies between the \textit{elements of a subset} and the
\textit{distinctions of a partition}, where a \textit{distinction} or
\textit{dit} of a partition $\pi=\left\{  B\right\}  $ on $U$ is an ordered
pair $\left(  u,u^{\prime}\right)  \in U\times U$ of elements in different
blocks of $\pi$. There is a set-theoretic representation of the lattice of
partitions where each partition is represented by its \textit{set of
distinctions} or \textit{ditset:}

\begin{center}
$\operatorname*{dit}\left(  \pi\right)  =\left\{  \left(  u,u^{\prime}\right)
:\exists B,B^{\prime}\in\pi;B\neq B^{\prime};u\in B;u^{\prime}\in B^{\prime
}\right\}  $,
\end{center}

\noindent and where the partial order is just inclusion between ditsets since: 

\begin{center}
$\sigma\preceq\pi$ iff $\operatorname*{dit}\left(  \sigma\right)
\subseteq\operatorname*{dit}\left(  \pi\right)  $. 
\end{center}

The complement of a ditset is the set of \textit{indistinctions} of the partition:

\begin{center}
$\operatorname*{indit}\left(  \pi\right)  =U\times U-\operatorname*{dit}%
\left(  \pi\right)  =\left\{  \left(  u,u^{\prime}\right)  :\exists B\in
\pi;u,u^{\prime}\in B\right\}  $
\end{center}

\noindent which is simply the equivalence relation associated with the
partition. The ditsets of partitions on $U$ are thus the complements of
equivalence relations on $U$ and they might be called the \textit{partition
relations} on $U$. 

Given any subset $S\subseteq U\times U$, its
\textit{reflexive-symmetric-transitive or rst closure} $cl\left(  S\right)  $
is the smallest equivalence relation containing $S$ (which is well-defined
since the intersection of two equivalence relations is an equivalence
relation). But it might be noted that this closure operation is not a
topological closure operation since the union of two rst-closed sets is not
necessarily rst-closed. The \textit{interior} $\operatorname*{int}\left(
S\right)  $ of a subset $S\subseteq U\times U$ is the complement of the
rst-closure of the complement, i.e., $\operatorname*{int}\left(  S\right)
=cl\left(  S^{c}\right)  ^{c}$, so it is the ditset of some partition. To
define the partition operation corresponding to any logical subset operation
(e.g., union, intersection, conditional, etc.), apply the subset operation to
the ditsets of the partitions, take the interior of the result, and then the
partition corresponding to that interior. For instance, the meet $\pi
\wedge\sigma$ of two partitions $\pi$ and $\sigma$ may be defined by the ditset:

\begin{center}
$\operatorname*{dit}\left(  \pi\wedge\sigma\right)  =\operatorname*{int}%
\left[  \operatorname*{dit}\left(  \pi\right)  \cap\operatorname*{dit}\left(
\sigma\right)  \right]  $.
\end{center}

Thus we can take any formula of subset logic and interpret it as a formula of
partition logic. The atomic variables would represent partitions on $U$
instead of subsets of $U$. Given such an interpretation of a formula
$\Phi\left(  \pi,\sigma,...\right)  $, an member $u\in U$ being an element of
the subset represented by $\Phi\left(  \pi,\sigma,...\right)  $ is analogous
to an ordered pair $\left(  u,u^{\prime}\right)  $ being a distinction of the
partition represented by $\Phi\left(  \pi,\sigma,...\right)  $. The two
definitions of a valid formula are also analogous. A formula $\Phi\left(
\pi,\sigma,...\right)  $ is a valid formula of subset logic, i.e., a
tautology, if for any subsets of $U$ substituted for the variables, the
formula evaluates to the set of all possible elements $U$ (the top of the
lattice) for any $U$ ($\left\vert U\right\vert \geq1$). Similarly, a formula
$\Phi\left(  \pi,\sigma,...\right)  $ is a \textit{valid formula of partition
logic}, i.e., a \textit{partition tautology}, if for any partitions on $U$
substituted for the variables, the formula evaluates to the partition that
makes all possible distinctions, i.e., the top-of-the-lattice discrete
partition $\mathbf{1}$ with the ditset $\operatorname*{dit}\left(  1\right)
=U\times U-\Delta$ ($\Delta$ is the diagonal $\left\{  \left(  u,u\right)
:u\in U\right\}  $), for any $U$ ($\left\vert U\right\vert \geq2$).

The following Figure $2$ summarizes the dual relationships between the two
logics (see \cite{ell:partitions} or \cite{ell:intropartlogic} for more on
partition logic).%

\begin{center}
\includegraphics[
height=2.7671in,
width=5.0892in
]%
{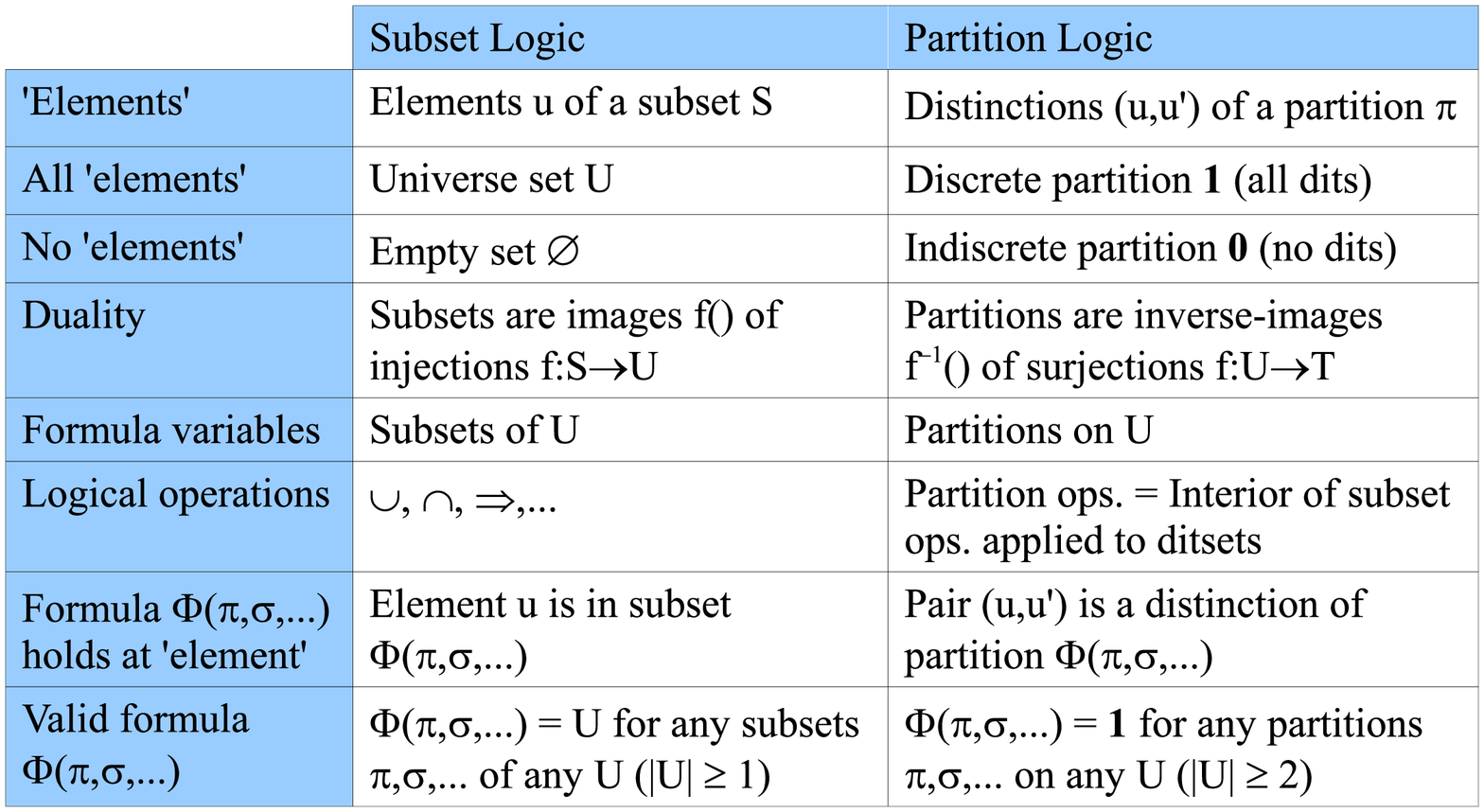}%
\\
Figure $2$: Dual analogies between the subset and partition logics
\end{center}

\section{Four ways to go from universal to particular}

In the two lattices of the dual logics, there are four universals, the tops
and bottoms of the two lattices. The four logico-mathematical ways to
abstractly characterize going from universal to particular are the four ways
of going from one of the universals to a particular subset or partition in the
lattice. In view of the duality between elements and distinctions, the four
ways can be characterized as: $\left(  1\right)  $ killing off elements,
$\left(  2\right)  $ creating elements, $\left(  3\right)  $ killing off
distinctions, and $\left(  4\right)  $ creating distinctions.

If $S$ represents any particular subset of $U$ and $\pi$ represents any
particular partition on $U$, then the four ways are:

\begin{enumerate}
\item the \textit{selectionist} mechanism $U\rightarrow S$: to go from the
universe set $U$ to a particular subset $S$ by "selecting" the elements of $S$
and eliminating or "killing off" the elements of the complement $S^{c}$;

\item the \textit{creationist} mechanism $\emptyset\rightarrow S$: to go from
the empty set $\emptyset$ to a particular subset $S$ by "creating" the
elements of $S$;

\item the \textit{identification} (or \textit{classification }or
\textit{symmetry-making}) mechanism $\mathbf{1}\rightarrow\pi$: to go from the
discrete partition $\mathbf{1}$ to a particular partition $\pi$ by identifying
elements ("killing off" distinctions) of $U$ (in a consistent way); and

\item the \textit{generative} (or \textit{symmetry-breaking}) mechanism
$\mathbf{0}\rightarrow\pi$: to go from the indiscrete partition $\mathbf{0}$
to a particular partition $\pi$ by "generating" distinctions on $U$ (in a
consistent way).
\end{enumerate}

The four schemes can be related in terms of duals (the elements-distinctions
duality) and opposites (all versus none) as in Figure $3$. %

\begin{center}
\includegraphics[
height=2.2208in,
width=3.9535in
]%
{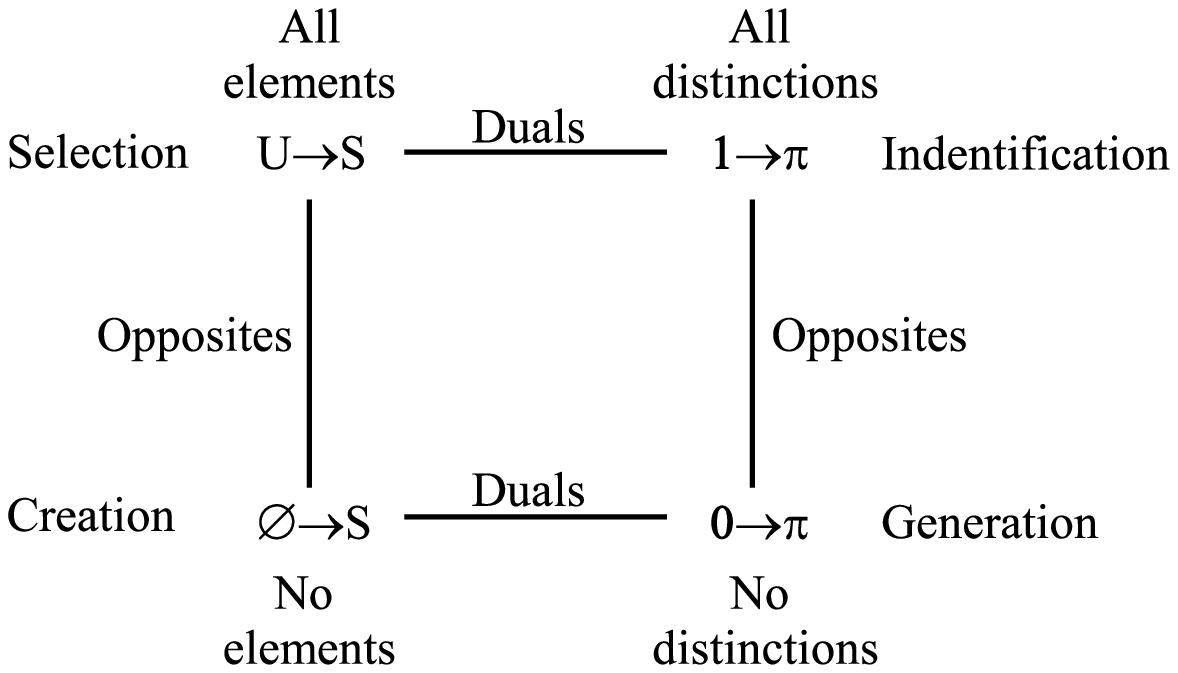}%
\end{center}

\begin{center}
Figure $3$: The dual \& opposite relations between the four schemes
\end{center}

\noindent A selectionist ($U\rightarrow S$) mechanism and a generative
($\mathbf{0}\rightarrow\pi$) mechanism are related by taking the dual
\textit{and} the opposite (in either order).

\subsection{The selectionist mechanism: From universe set to a subset}

The first scheme $U\rightarrow S$ is the abstract logico-mathematical model of
a selectionist process since it starts with an actualized set of diverse
alternatives $U$ and then a number of the alternatives are eliminated by some
fitness criterion or evaluation metric while the remaining alternatives are
selected (e.g., by differential amplification).%

\begin{center}
\includegraphics[
height=1.8962in,
width=1.863in
]%
{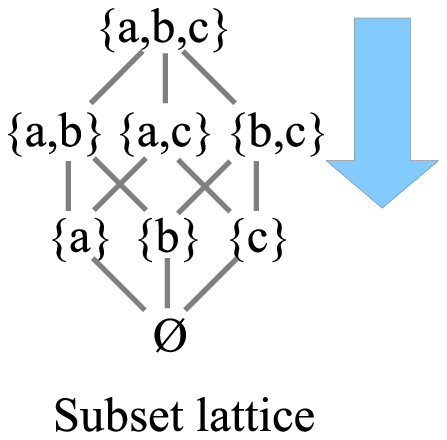}%
\end{center}

\begin{center}
Figure $4$: From the universe set $U$ to a particular subset $S$
\end{center}

The original example of a selectionist process is Darwinian evolution where
the set of diverse alternatives is generated over time by random genetic
mutations and then the environment applies a fitness filter (see
\cite{fodorp-p:darwinwrong} on the general generate-and-filter idea). In the
application to the immune system, the actualized set of diverse alternatives
are the generated set of antibodies in low concentration and then the
selectionist process differentially amplifies those that fit an invading
antigen. Gerald Edelman's various selectionist models \cite{edel:nd} also fit
this scheme. One of Peter Medawar's metaphors for a selectionist scheme was a
jukebox where all the tunes were already actualized as records in the jukebox
and then one was selected. In the format-selection or FS approach to universal
grammar, the mental mechanism must generate some representation of the diverse
variety of grammars, and then a chunk of linguistic experience is evaluated
according to some evaluation metric to find the best fit among the various
systems of rules.

\subsection{The creationist mechanism: From empty set to a subset}

The second scheme $\emptyset\rightarrow S$ is the abstract logico-mathematical
model of a creation story where elements are, in effect, created out of
nothing. While this is one type of "creation story," the Big Bang creation
theory \cite{pagels:perfsymm} is modelled not by this $\emptyset\rightarrow S$
scheme but by the generative $\mathbf{0}\rightarrow\pi$ scheme where the
making of distinctions is rendered as symmetry-breaking. The $\emptyset
\rightarrow S$ scheme is perhaps the least interesting to model actual processes.%

\begin{center}
\includegraphics[
height=1.8962in,
width=1.863in
]%
{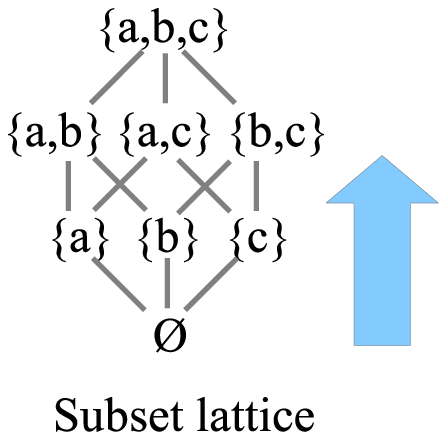}%
\end{center}

\begin{center}
Figure $5$: From the empty set $\emptyset$ to a particular subset $S$
\end{center}

\subsection{The identification mechanism: From discrete partition to a
partition}

The third scheme $\mathbf{1}\rightarrow\pi$ is the abstract
logico-mathematical model of any classification, partitioning \cite[p.
82]{law:conceptual}, or quotienting process that proceeds by making consistent
identifications ("consistent" means that the identifications must be
reflexive, symmetric, and transitive to form an equivalence relation).%

\begin{center}
\includegraphics[
height=1.8032in,
width=2.6094in
]%
{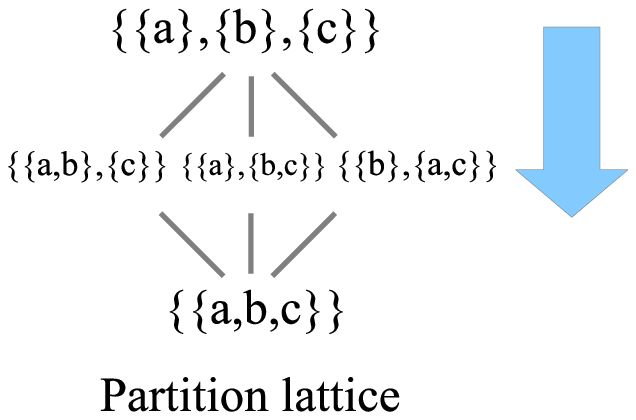}%
\end{center}

\begin{center}
Figure $6$: From the discrete partition $\mathbf{1}$ to a particular partition
$\pi$
\end{center}

For instance, the classification of animals to species where $\mathbf{1}$
represents each animal by itself and $\pi$ represents the partition of the set
of animals as to species. Mathematically, the action of a group on a set is
automatically reflexive, symmetric, and transitive so it defines an
equivalence relation where the equivalence classes are called "orbits"
\cite[p. 99]{macbirk:algebra}. This $\mathbf{1}\rightarrow\pi$ scheme is
"symmetry-making" while the opposite scheme $\mathbf{0}\rightarrow\pi$ is "symmetry-breaking."

\subsection{The generative mechanism: From indiscrete partition to a
partition}

The fourth scheme $\mathbf{0}\rightarrow\pi$ is the abstract
logico-mathematical model of any \textit{generative} process where a number of
different outcomes (represented by the blocks of $\pi$) can be generated by
consistently making distinctions ("consistent" means nothing can be
distinguished from itself, distinguishing must be symmetric, and if $u$ is
distinguished from $u^{\prime}$ and $u=u_{1},u_{2},...,u_{n}=u^{\prime}$, then
one of the pairs $\left(  u_{i},u_{i+1}\right)  $ must also be distinguished
for $i=1,...,n-1$, all of which means the set of distinctions must be
anti-reflexive, symmetric, and anti-transitive, i.e., a partition relation).%

\begin{center}
\includegraphics[
height=1.8032in,
width=2.6094in
]%
{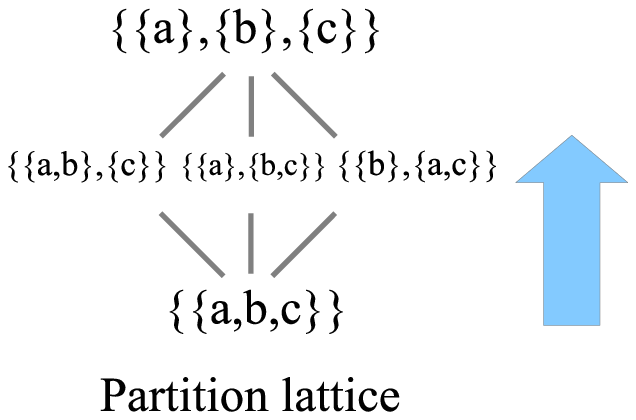}%
\end{center}

\begin{center}
Figure $7$: From the indiscrete partition $\mathbf{0}$ to a particular
partition $\pi$
\end{center}

The most "atomic" type of distinction is a binary partition of the single
block $\left\{  U\right\}  \in\mathbf{0}$ into two blocks and then the binary
partitions can be combined or joined together. The \textit{join} $\pi
\vee\sigma$ of two partitions $\pi=\left\{  B\right\}  $ and $\sigma=\left\{
C\right\}  $ is the partition of nonempty intersections $B\cap C$ (or, in
terms of ditsets, $\operatorname*{dit}\left(  \pi\vee\sigma\right)
=\operatorname*{dit}\left(  \pi\right)  \cup\operatorname*{dit}\left(
\sigma\right)  $, where the interior is not needed since a union of partitions
relations is always a partition relation). And the most "efficient" binary
partition is one that divides the block $\left\{  U\right\}  $ into two equal
parts (assuming an even number of elements). The classic example is where $U$
has $2^{n}$ elements which can be enumerated using $n$-place binary numbers.
Then $U$ can be divided into two equal parts by the binary partition according
to whether the $i^{th}$ binary digit is $0$ or $1$. The join of those binary
partitions for $i=1,...,n$ would go from the indiscrete partition $\mathbf{0}$
all the way to the discrete partition $\mathbf{1}$, so the $n$ equal-binary
partitions are Shannon's $n$ \textit{bits} \cite{ell:distinctions}.

Often a $\mathbf{0}\rightarrow\pi$ generative process proceeds not only by
joining binary partitions (with not necessarily equal blocks) but by
designating one of the blocks as in the game of twenty questions where the
block with the yes-answer to the yes-or-no question is designated. In this
case, the $\mathbf{0}\rightarrow\pi$ process goes not just from the indiscrete
partition $\mathbf{0}$ to a particular partition $\pi$ but from the single
block $\left\{  U\right\}  \in\mathbf{0}$ to a specific block $B\in\pi$ (like
a correct answer in the game of twenty questions) by following the
yes-branches on the binary tree. This is the case of the generative mechanism
that is of most interest for our purposes.

A binary partition with a designated block is just a choice with two options,
and it might be represented by a switch with a neutral setting (representing
the state before the choice is made) and then two options such as a Left
Option and a Right Option.%

\begin{center}
\includegraphics[
height=1.5235in,
width=4.8817in
]%
{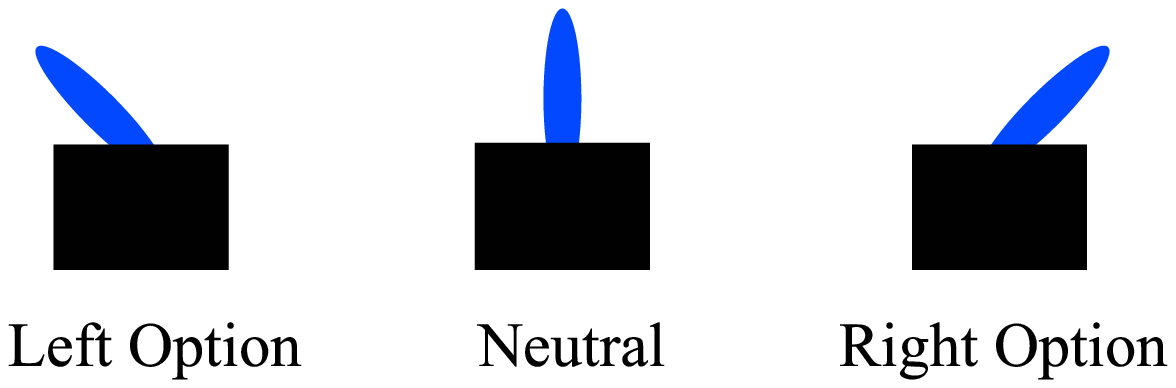}%
\end{center}

\begin{center}
Figure $8$: Switch to go from Neutral to either Left or Right Options
\end{center}

For instance consider the example where $U$ has $8=2^{3}$ elements represented
by the $3$-digit binary numbers $b_{3}b_{2}b_{1}$. There would be three
switches for $i=1,2,3$ where the two options for the $i^{th}$ switch are a $0$
or $1$ in the $i^{th}$ place $b_{i}$ of the $3$-digit binary number. Then the
initial state is the indiscrete partition $\mathbf{0}=\left\{  U\right\}  $
where all the switches are in neutral. The first switch (controlling the first
digit on the right) determines the binary partition with one block having the
four elements $b_{3}b_{2}0$ (Left Option) and the other block having the four
elements $b_{3}b_{2}1$ (Right Option), and so forth for the other two
switches. When all three switches have been set one way or the other, that
determines the transition from the single block $\left\{  U\right\}
\in\mathbf{0}$ to a specific singleton block represented by a specific
$3$-digit binary number.

Our principal application of the $\mathbf{0}\rightarrow\pi$ generative scheme
is Chomsky's principles and parameters (P\&P) description of the
language-acquisition faculty or universal grammar (\cite{chom:govt},
\cite{chomlas:pandp}, \cite{chom:minprog}). In our simple model, the
parameters are represented by the switches that can be moved to the left or
right (from the original setting of neutral) by the child's linguistic
experience, and the underlying principles are expressed in the whole setup
defining the grammatical meaning of the left and right settings.

\begin{quotation}
A simple image may help to convey how such a theory might work. Imagine that a
grammar is selected (apart from the meanings of individual words) by setting a
small number of switches--0, say--either "On" or "Off." Linguistic information
available to the child determines how these switches are to be set. In that
case, a huge number of different grammars (here, 2 to the twentieth power)
will be prelinguistically available, although a small amount of experience may
suffice to fix one. \cite[p. 154]{higg:chomsky}
\end{quotation}

\noindent Needless to say, this imagery implicitly allows for a neutral
setting on the switches (sometimes called the "initial state $S_{0}$") since
otherwise the original or "factory" setting of the switches would determine a
specific grammar independent of experience.

\begin{quotation}
\noindent Borrowing an image suggested by James Higginbotham, we may think of
UG as an intricately structured system, but one that is only partially "wired
up."\ The system is associated with a finite set of switches, each of which
has a finite number of positions (perhaps two). Experience is required to set
the switches. When they are set, the system functions. The transition from the
initial state $S_{0}$ to the steady state $S_{s}$ is a matter of setting the
switches. \cite[p. 146]{chom:knowoflang}
\end{quotation}

Another implication of this general type of $\mathbf{0}\rightarrow\pi$ model
for Chomsky's language-acquisition faculty is the interpretation of the
adjective "universal" in the phrase "universal grammar." It does \textit{not}
mean a specific grammatical rule common to all languages (which would be a
subset-logic interpretation of "universal").

\begin{quotation}
\noindent The switch-settings of the metaphor above are in Chomsky's
terminology the "parameters" defined by universal grammar. Notice that this
image underscores the sense in which universal grammar, the initial state of
the language-learner, need not comprise an account of what languages have in
common--to continue the metaphor, different switch-settings could give rise to
very different grammatical system. \cite[p. 154]{higg:chomsky}
\end{quotation}

\section{Generative versus selectionist mechanisms}

It might also be useful to illustrate a selectionist and a generative
mechanism to solve the same problem of determining one among the $8=2^{3}$
options considered in the last section. The eight possible outcomes might be
represented as:

\begin{center}
$\left\vert 000\right\rangle $, $\left\vert 100\right\rangle $, $\left\vert
010\right\rangle $, $\left\vert 110\right\rangle $, $\left\vert
001\right\rangle $, $\left\vert 101\right\rangle $, $\left\vert
011\right\rangle $, $\left\vert 111\right\rangle $.
\end{center}

\noindent In the selectionist scheme, all eight variants are in some sense
actualized or realized in the initial state $S_{0}$ so that a fitness
criterion or evaluation metric (as in the FS scheme) can operate on them. Some
variants do better and some worse as indicated by the type size in Figure $9$.%

\begin{center}
\includegraphics[
height=2.3694in,
width=4.5313in
]%
{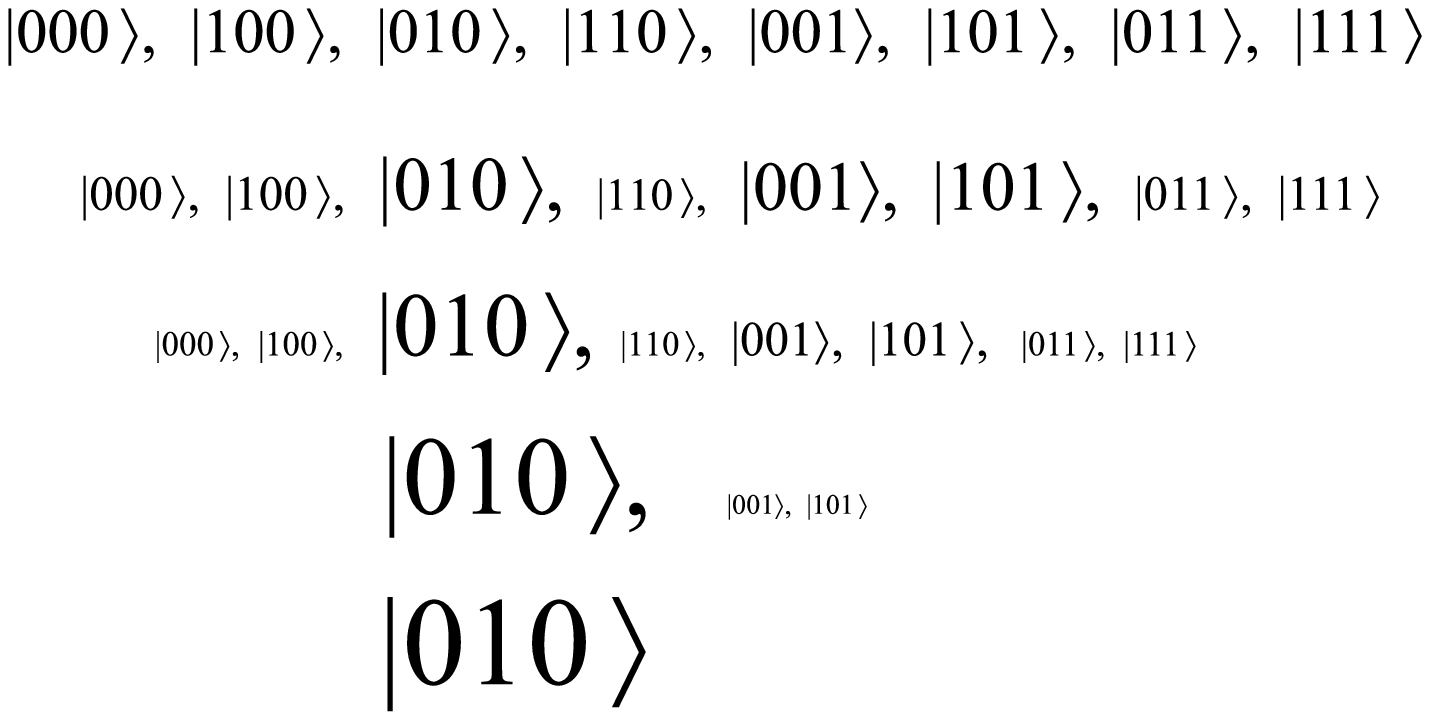}%
\end{center}

\begin{center}
Figure $9$: An abstract model of a selectionist learning mechanism
\end{center}

\noindent Eventually the "unfit" options dwindle, atrophy, or die off leaving
the most fit option $\left\vert 010\right\rangle $ as the final steady state
$S_{s}$.

In the generative learning scheme, the initial state $S_{0}$ is where all the
switches are in neutral so all the eight potential outcomes are in a
"superposition" state indicated by the plus signs in the following Figure $10$.%

\begin{center}
\includegraphics[
height=1.8107in,
width=4.6094in
]%
{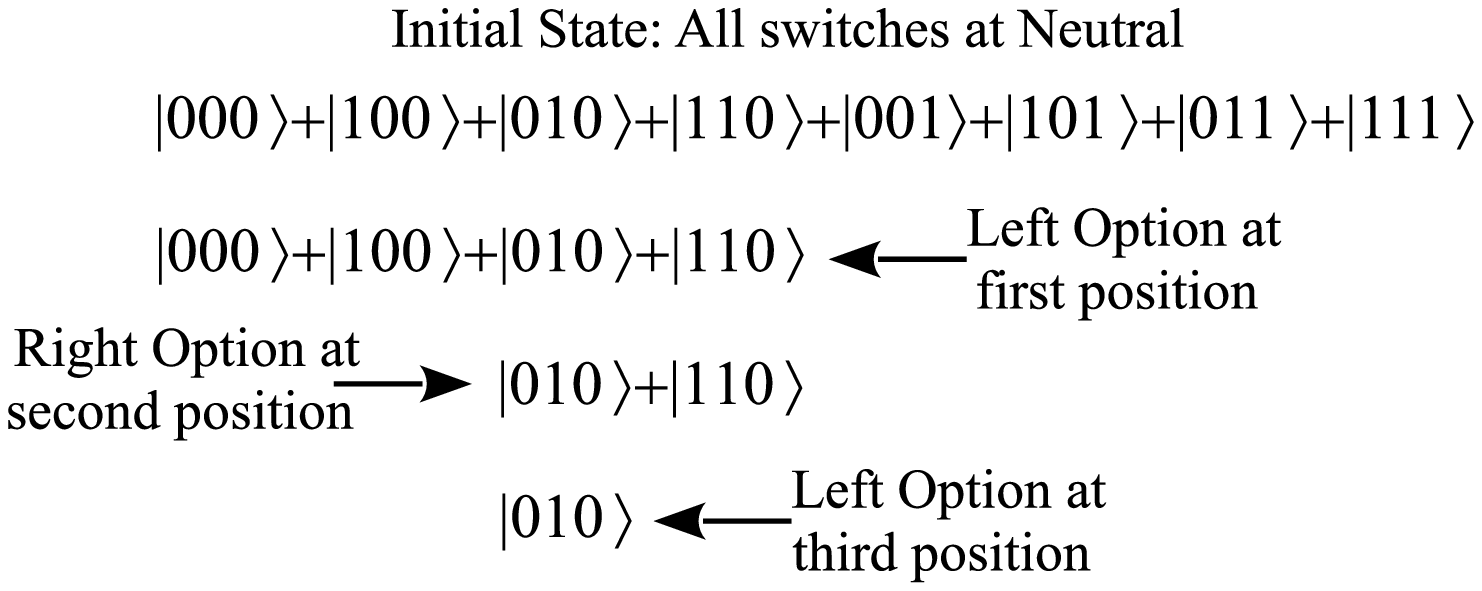}%
\end{center}

\begin{center}
Figure $10$: An abstract model of a generative learning mechanism
\end{center}

\noindent We assume that the initial experience sets the first switch (or
first parameter in the P\&P model) to the left option which reduces the state
to $\left\vert 000\right\rangle +\left\vert 100\right\rangle +\left\vert
010\right\rangle +\left\vert 110\right\rangle $ (where the plus signs between
these options indicate that the second and third switches are still in
neutral). Then subsequent experience sets the second switch to the right
option and the third switch to the left option. Thus we reach the same outcome
$\left\vert 010\right\rangle $ as the final state $S_{s}$ in the two models
but by quite different mechanisms.

We have now differentiated a selectionist ($U\rightarrow S$) mechanism from a
generative ($\mathbf{0}\rightarrow\pi$) mechanism at a very abstract
logico-mathematical level. 

There is a whole literature where "selectionist" is interpreted very broadly
as non-instructionist so that the generative ($\mathbf{0}\rightarrow\pi$)
mechanism is also described in those overly-broad terms as being
"selectionist" (e.g., by describing the generative mechanism as "selecting" switch-settings).

\begin{quotation}
I suggest that some important lessons for linguistics and cognitive science
can, indeed, be drawn from contemporary biology, but that the new principles
and the new assumptions came to bury learning by instruction and to replace It
with learning by selection, a radically different process. What now replaces
learning everywhere in biology has nothing to do with a transfer of structure
and everything to do with mechanisms of internal selection and filtering
affecting a pre-programmed chain of multiple internal recombinations and
internal "switches." \cite[p. 3]{piattell:selection}
\end{quotation}

Now we see that this sort of setting of internal switches is better described
as a generative ($\mathbf{0}\rightarrow\pi$) mechanism whereas many of the
other "learning" mechanisms in biology (e.g., in the immune system) are
correctly described as a selectionist ($U\rightarrow S$) mechanism since the
latter involves the actualization of some "universal" repertoire of
possibilities some of which are selected.

A hierarchy of genetic switches as in a stem cell or in embryonic development
would be a generative mechanism (\cite{monad:chance}, \cite{jacob:logic}). One
might imagine a hypothetical selectionist mechanism to replace stem cells that
would postulate low concentrations of the different types of cells through the
body, so that, say, muscle cells would be selected to multiply in a muscle
environment while the other types of cells would be inactive there. Yet what
is found biologically is not that type of selectionist mechanism but the
generative mechanism of stem cells (where the muscle environment sets the
switches to produce a muscle cell--in addition to reproducing the stem cell).

Peter Metawar \cite{med:fut} explains the selectionist-instructionist
juxtaposition by contrasting a jukebox (with the musical records taken as
internal) with a record player (with the records taken as external). The
jukebox has a set of pre-existing options one of which is selected by the
simple pushing of a button whereas when a record player plays music, the set
of external instructions must be supplied in the form of a record. Thus a
jukebox is a genuine selectionist ($U\rightarrow S$) mechanism. Medawar also
describes the development of the embryo as being selectionist:

\begin{quotation}
\noindent Embryonic development\ldots\ must therefore be an unfolding of
pre-existing capabilities, an acting-out of genetically encoded instructions;
the inductive stimulus is the agent that selects or activates one set of
instructions rather than another. \cite[p. 295]{med:plato}
\end{quotation}

\noindent But in terms of our differentiation, embryonic development is a
generative ($\mathbf{0}\rightarrow\pi$), \textit{not} a selectionist
($U\rightarrow S$), mechanism. In a similar manner, it is easy to see that
Gerald Edelman's \cite{edel:nd} various models of brain development and
learning are all selectionist ($U\rightarrow S$) mechanisms.

In this manner, one could go over all the examples broadly called
"selectionist" and see which were genuinely selectionist ($U\rightarrow S$)
mechanisms and which were generative ($\mathbf{0}\rightarrow\pi$)
mechanisms--which shows the surprising fruitfulness of the quite abstract
logico-mathematical differentiation between $U\rightarrow S$ and
$\mathbf{0}\rightarrow\pi$ mechanisms.

\end{document}